\documentclass[11pt]{article}
\makeindex
\usepackage{subfigure}

\usepackage{
amssymb,amsmath,graphics,latexsym,amsfonts,indentfirst}
\usepackage{epsf,epsfig,color, pdfpages}
\usepackage{amssymb,latexsym}
\usepackage{graphicx}%
\usepackage{amsmath,bm}
\usepackage{graphicx}
\usepackage{color}
\usepackage{ifpdf}
\usepackage{epstopdf}
\usepackage{pinlabel}
\usepackage[noadjust]{cite}
\usepackage{etoolbox}
\usepackage{mathrsfs}
\usepackage{lipsum}

\newtheorem{lem}{Lemma}
\newtheorem{theo}{Theorem}
\newtheorem{pro}{Proposition}
\newtheorem{con}{Conjecture}
\newtheorem{cor}{Corollary}

\newcommand \ARRAY[1]
{
\left \{
\begin{array}{ll}
#1
\end{array}
\right.
}

\newcommand \lemma[2]
{\begin{lem}
\label{#1} #2
\end{lem}
}

\newcommand \prop[2]
{\begin{pro}
\label{#1} #2
\end{pro}
}

\newcommand{\proof}{{\noindent {\em Proof}.\quad}\setcounter{countclaim}{0}\setcounter{countcase}{0}}

\newcommand{\proofend}{{\hfill$\Box$}}

\newcounter{countcase}

\newcounter{countclaim}
\def\inclaim{\addtocounter{countclaim}{1}
{\noindent {\bf Claim \thecountclaim}: }}

\newcounter{fourclaim-1}
\newcounter{fourclaim-2}

\newcounter{countfig}

\newcommand{\beeq}{\begin{equation}}
\newcommand{\eneq}{\end{equation}}

\newcommand{\beeqn}{\begin{eqnarray*}}
\newcommand{\eneqn}{\end{eqnarray*}}

\def \T {\mathbb{T}}

\def \SEDF {{\cal F}_{sed}}

\begin{document}

\baselineskip 0.55 cm

\title
{Upper bounds on the signed edge domination number of a graph}

\author{
Fengming Dong$^{1}$\thanks{Email: fengming.dong@nie.edu.sg and donggraph@163.com.},
Jun Ge$^{2}$\thanks{Corresponding author. Email: mathsgejun@163.com.},
and Yan Yang$^{3}$\thanks{Email: yanyang@tju.edu.cn.}\\[2mm]
\small $^{1}$National Institute of Education, Nanyang Technological University, Singapore \\
\small $^{2}$School of Mathematical Sciences, Sichuan Normal University, China \\
\small $^{3}$School of Mathematics, Tianjin University, China.
}

\date{}

\maketitle{}

\begin{abstract}
A signed edge domination function (or SEDF) of a simple graph
$G=(V,E)$ is a function $f: E\rightarrow \{1,-1\}$
such that $\sum_{e'\in N[e]}f(e')\ge 1$ holds for each
edge $e\in E$, where $N[e]$ is the set of edges in $G$
that share at least one endpoint with $e$.
Let $\gamma_s'(G)$ denote the minimum value of $f(G)$
among all SEDFs $f$, where $f(G)=\sum_{e\in E}f(e)$.
In 2005, Xu conjectured that $\gamma_s'(G)\le n-1$, where $n$ is the order of $G$.
This conjecture has been proved
for the two cases $v_{odd}(G)=0$ and $v_{even}(G)=0$,
where $v_{odd}(G)$ (resp. $v_{even}(G)$) is the number of odd (resp. even) vertices in $G$.
This article proves
Xu's conjecture for $v_{even}(G)\in  \{1, 2\}$.
We also show that
for any simple graph $G$ of order $n$,
$\gamma_s'(G)\le n+v_{odd}(G)/2$ and
$\gamma_s'(G)\le n-2+v_{even}(G)$ when $v_{even}(G)>0$,
and thus $\gamma_s'(G)\le (4n-2)/3$.
Our result improves the best current upper bound of
$\gamma_s'(G)\le \lceil 3n/2\rceil$.
\end{abstract}

\medskip



\section{Introduction}

This article considers simple and undirected graphs only.
For a graph $G$, let $V(G)$ and $E(G)$ denote
its vertex set and edge set, respectively.
For any $v\in V(G)$,
let $E_G(v)$ be the set of edges in $G$ incident to $v$,
let $N_G(v)$  be the set
of vertices in $G$ adjacent to $v$,
and let $N_G[v]=N_G(v)\cup \{v\}$.
$E_G(v)$, $N_G(v)$ and $N_G[v]$
are simply written as $E(v)$, $N(v)$ and $N[v]$,
respectively, when there is no confusion.
For any $v\in V(G)$, we use $d_G(v)$ (or simply $d(v)$ when there is no confusion)
to denote the degree of $v$ in $G$.

For a graph $G=(V,E)$,
a {\it signed domination function} of $G$ is a function
$f:V\rightarrow \{1,-1\}$ with the property that
$f(N[v])\ge 1$ holds for every $v\in V$,
where $f(S)=\sum_{v\in S} f(v)$ for each $S\subseteq V$.
The {\it signed domination number} of $G$,
denoted by $\gamma_s(G)$,
is defined to be the minimum value of $f(V)$
over all signed domination functions $f$ of $G$.
The parameter $\gamma_s(G)$ was introduced by
Dunbar, Hedetniemi, Henning, and
Slater~\cite{dunbar1995signed}
and has been studied by many authors,
e.g., \cite{chen2008lower,favaron1996signed,furedi1999signed,henning1999open,matouvsek2000signed,zelinka1996some}.

In 2001, Xu \cite{xu2001signed} introduced
signed edge domination functions.
For a graph $G=(V,E)$,
a function $f :E\rightarrow \{1,-1\}$
is called a {\it signed edge domination function} (SEDF)
of $G$ if $\sum_{e'\in N[e]} f(e')\ge 1$ holds
for every $e\in E$,
where $e=uv$ and $N[e]=E_G(u)\cup E_G(v)$.
Let $\SEDF(G)$ denote the set of SEDFs of $G$.
The {\it signed edge domination number} of $G$, denoted by
$\gamma'_s(G)$, is defined to be the minimum value of $f(G)$
over all $f\in \SEDF(G)$, where $f(G)=\sum_{e\in E}f(e)$.

Observe that the parameter $\gamma'_s(G)$ is an extension of
$\gamma_s(G)$,
as each member $f$ in $\SEDF(G)$
is actually a signed domination function of the line graph $L(G)$,
thus implying that $\gamma'_s(G)=\gamma_s(L(G))$.
The parameter $\gamma'_s(G)$
has been studied by many authors,
e.g.,
\cite{akbari2009signed,akbari2014some,karami2009improved,xu2001signed,xu2005edge,xu2006two,xu2009signed}.
The following are some known results on $\gamma'_s(G)$
for a graph $G$ of order $n$ and size $m$:
\begin{enumerate}

\item $\gamma'_s(G)\ge \frac{-n^2}{16}$ \cite{akbari2009signed};

\item for any positive integer $r$, there exists
an $r$-connected graph $H$ such that
$\gamma'_s(H)\le -\frac {r}6 |V(H)|$ \cite{akbari2009signed};

\item $\gamma'_s(G)\ge \frac{2\alpha'(G)-m}{3}$,
where $\alpha'(G)$ is the size of
a largest matching of $G$ \cite{akbari2014some};

\item $\gamma'_s(G)\ge n-m$ for $n\ge 4$ \cite{xu2005edge};

\item $\gamma'_s(G)\le \frac{11n}6-1$ \cite{xu2006two};

\item $\gamma'_s(G)\le \lceil \frac{3n}2\rceil$ \cite{karami2009improved}.
\end{enumerate}

In this article, we will improve the upper bounds of $\gamma_s'(G)$
by establishing the following result.
A vertex in a graph $G$ is called an
{\it odd vertex} (resp. {\it even vertex}) if it is
of odd degree (resp. even degree) in $G$.
Let $v_{odd}(G)$ (resp. $v_{even}(G)$) denote the
number of odd (resp. even) vertices in $G$.
Clearly, $v_{odd}(G)$ is even.

\begin{theo}\label{th2}
For any graph $G$ of order $n$,

(a) $\gamma_s'(G)\le n+v_{odd}(G)/2$;

(b) $\gamma_s'(G)\le n-2+v_{even}(G)$ when $v_{even}(G)>0$;

and hence $\gamma_s'(G)\le (4n-2)/3$.
\end{theo}

The most challenging and interesting problem on
$\gamma'_s(G)$ may be the following conjecture proposed by
Xu \cite{xu2005edge} in 2005.

\begin{con}[\cite{xu2005edge}]\label{conj1}
For any simple graph $G$ of order $n$,
$\gamma_s'(G)\le n-1$ holds.
\end{con}

As far as we know, Conjecture~\ref{conj1} has been
only proved for a few cases.
Karami, Khodkar, and Sheikholeslami~\cite{karami2009improved} showed
that Conjecture~\ref{conj1} holds
when $v_{odd}(G)\in \{0,n\}$.
In the case  $v_{odd}(G)=n$,
Akbari, Esfandiari, Barzegary, and Seddighin~\cite{akbari2014some}
strengthened the result to
$\gamma_s'(G)\le n-\frac{2\alpha'(G)}{3}$,
where $\alpha'(G)$ is the size of a maximum matching in $G$.
In this paper, we prove that $\gamma_s'(G)\leq n-1$
if $v_{even}(G)\in \{1, 2\}$.

\begin{theo}\label{th3}
Conjecture $\ref{conj1}$ holds for any simple graph
$G$ with $v_{even}(G)\in \{1, 2\}$.
\end{theo}

In Section~\ref{sec2}, we
introduce a subfamily $\SEDF^0(G)$  of
$\SEDF(G)$ and establish some basic results
for proving the main results in
the following sections.
Theorem~\ref{th2} (a) and (b)
are proved in Sections~\ref{sec3} and~\ref{sec4},
respectively.
By Theorem~\ref{th2} (b),
Conjecture~\ref{conj1} holds for $v_{even}(G)=1$.
In Section~\ref{sec5}, we show that
Conjecture~\ref{conj1} holds for $v_{even}(G)=2$,
and thus Theorem~\ref{th3} follows.
In Section~\ref{sec6},
we propose a conjecture to replace
Conjecture~\ref{conj1}, as we think
there exists
a member $f$ in
$\SEDF^0(G)$ with $f(G)\le n-1$ for any
graph $G$ of order $n$.
We also propose a conjecture for
the lower bound of $\gamma_s'(G)$ when $G$ is $2$-connected.

\section{A subset $\SEDF^0(G)$ of $\SEDF(G)$
\label{sec2}
}

Let $G$ be a simple graph.
For any $f:E(G)\rightarrow \{1,-1\}$ and $v\in V(G)$,
let $f(v)=\sum\limits_{e\in E_G(v)}f(e)$
and let $f(S)=\sum\limits_{e\in S}f(e)$, where $S\subseteq E(G)$.
Let $\SEDF^0(G)$ denote the set of functions
$f:E(G)\rightarrow \{1,-1\}$ satisfying the two conditions below:

(a) $f(v)\ge 0$ for all $v\in V(G)$; and

(b) $f(u)+f(v)\ge 2$ for each $e=uv\in E(G)$ with $f(e)=1$.

\lemma{le2-1}
{$\SEDF^0(G)\subseteq \SEDF(G)$.}

\proof Let $f$ be any member in $\SEDF^0(G)$ and
let $e=v_1v_2\in E(G)$.
It follows from the definition of $\SEDF^0(G)$
that $f(v_i)\ge 0$ for $i=1,2$ and
$f(v_1)+f(v_2)\ge 2$ holds whenever $f(e)=1$,
thus implying $f(v_1)+f(v_2)\ge 1+f(e)$.
Consequently,
$$
f(N[e])=f(v_1)+f(v_2)-f(e)\ge 1+f(e)-f(e)=1.
$$
Hence $f\in \SEDF(G)$ and the result holds.
\proofend

For $S\subseteq V(G)$,
let $G[S]$ denote the subgraph of $G$ induced by $S$.
If $E_1$ and $E_2$ form a partition of $E(G)$
and $f_i:E_i\rightarrow \{1,-1\}$,
let $f_1*f_2$ be the function $f:E(G)\rightarrow \{1,-1\}$
defined by
$f(e)=f_i(e)$ whenever $e\in E_i$.

\lemma{le2-2}
{Let $G$ be a separable graph with $V(G)=V_1\cup V_2$,
$V_1\cap V_2=\{v_0\}$
and $E(G)=E(G[V_1])\cup E(G[V_2])$.
If $f_i\in \SEDF^0(G[V_i])$ for $i=1,2$,
then $f=f_1*f_2\in \SEDF^0(G)$
with
$$
f(G)=f_1(G[V_1])+f_2(G[V_2]).
$$
}

\proof Note that $E(G[V_1])$ and $E(G[V_2])$ form
a partition of $E(G)$ and thus $f_1*f_2$ is well defined.
By the definition of $f$, it is obvious that
$f(G)=f_1(G[V_1])+f_2(G[V_2])$.
Next, by the definition of $f$,
for any $v\in V(G)$,
\begin{equation}\label{eq2-1}
f(v)=
\ARRAY
{
f_1(v_0)+f_2(v_0),\qquad &\mbox{if }v=v_0;\\
f_i(v),  &\mbox{if }v\in V_i-\{v_0\}, i=1,2.\\
}
\end{equation}
As $f_i\in \SEDF^0(G[V_i])$ for $i=1,2$,
we have $f(v)\ge 0$ for each $v\in V(G)$
by (\ref{eq2-1}).
Now, let $e$ be any edge in $E(G)$ with $f(e)=1$.
We may assume that $e=v_1v_2\in E(G[V_1])$,
and thus $f_1(e)=f(e)=1$.
As  $f_1\in \SEDF^0(G[V_1])$, $f_1(v_1)+f_1(v_2)\ge 2$.
By (\ref{eq2-1}) and the assumption that
$f_2\in \SEDF^0(G[V_2])$, we have
$$
f(v_1)+f(v_2)\ge f_1(v_1)+f_1(v_2)\ge 2.
$$
Hence $f\in \SEDF^0(G)$ as required.
\proofend

In the following, we assume that $v_0$ is
a vertex in a 2-connected graph $G$ with
$d_G(v_0)=2$.

\lemma{le2-3}
{Let $G$ be a simple graph, and let $v_0\in V(G)$
with $N_G(v_0)=\{u_1,u_2\}$.
\begin{enumerate}
\item For $u_1u_2\in E(G)$ and $g\in \SEDF^0(G')$, where $G'=G-u_1u_2-v_0$, as shown in Figure~\ref{f4}(b),
let $f: E(G)\rightarrow \{1,-1\}$ be defined below:
$$
f(e)=
\ARRAY
{g(e), \quad &\mbox{if }e\in E(G');\\
1,  &\mbox{if }e=v_0u_i,~i=1,2;\\
-1, &\mbox{if } e=u_1u_2.
}
$$
Then, $f\in \SEDF^0(G)$ with $f(G)=g(G')+1$.

\begin{figure}[ht!]
\centering

\includegraphics[width=12cm]{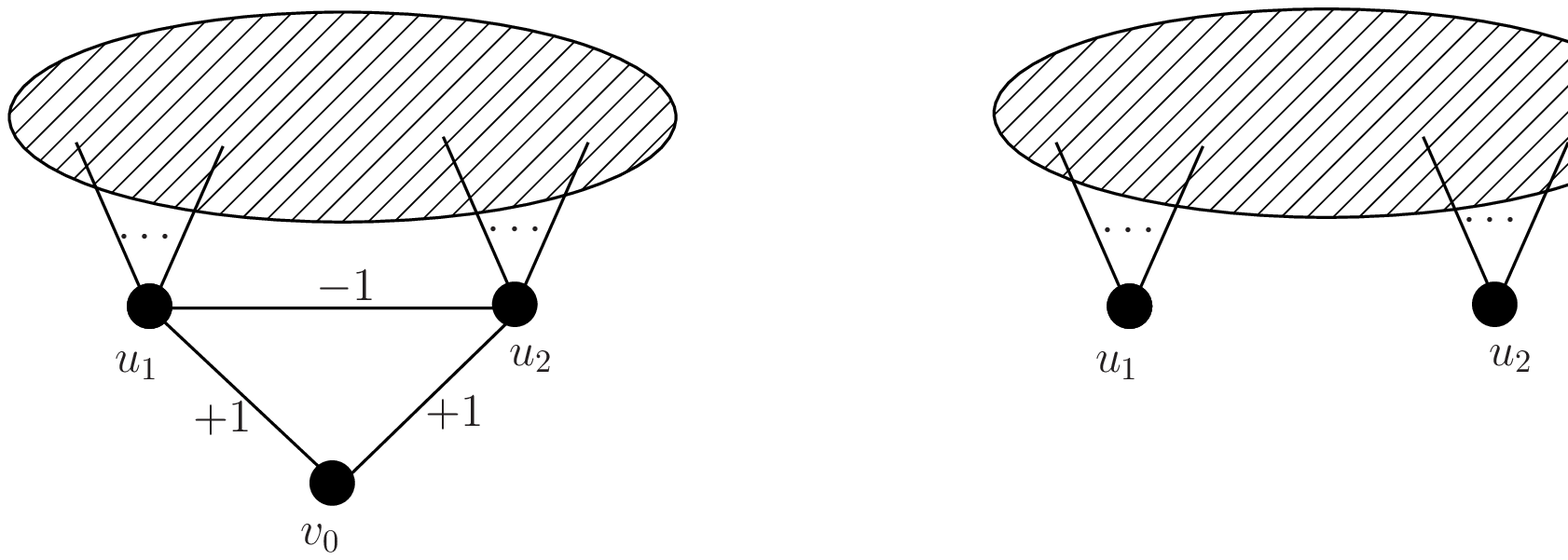}\\

\hspace{1 cm}(a) $G$ \hspace{4.5 cm} (b) $G'$ ($=G-u_1u_2-v_0$)

\caption{Graphs $G$ and $G'$ ($=G-u_1u_2-v_0$).}
\label{f4}
\end{figure}

\item For $u_1u_2\notin E(G)$ and $g\in \SEDF^0(G')$,
where $G'=G+u_1u_2-v_0$,
as shown in Figure~\ref{f5}(b),
let $f: E(G)\rightarrow \{1,-1\}$ be defined below:
$$
f(e)=
\ARRAY
{g(e), \quad &\mbox{if }e\in (E(G')-\{u_1u_2\});\\
1,  &\mbox{if }e=u_1v_0;\\
g(u_1u_2), \quad &\mbox{if }e=u_2v_0.
}
$$
Then, $f\in \SEDF^0(G)$ with $f(G)=g(G')+1$.

\begin{figure}[ht!]
\centering

\includegraphics[width=12cm]{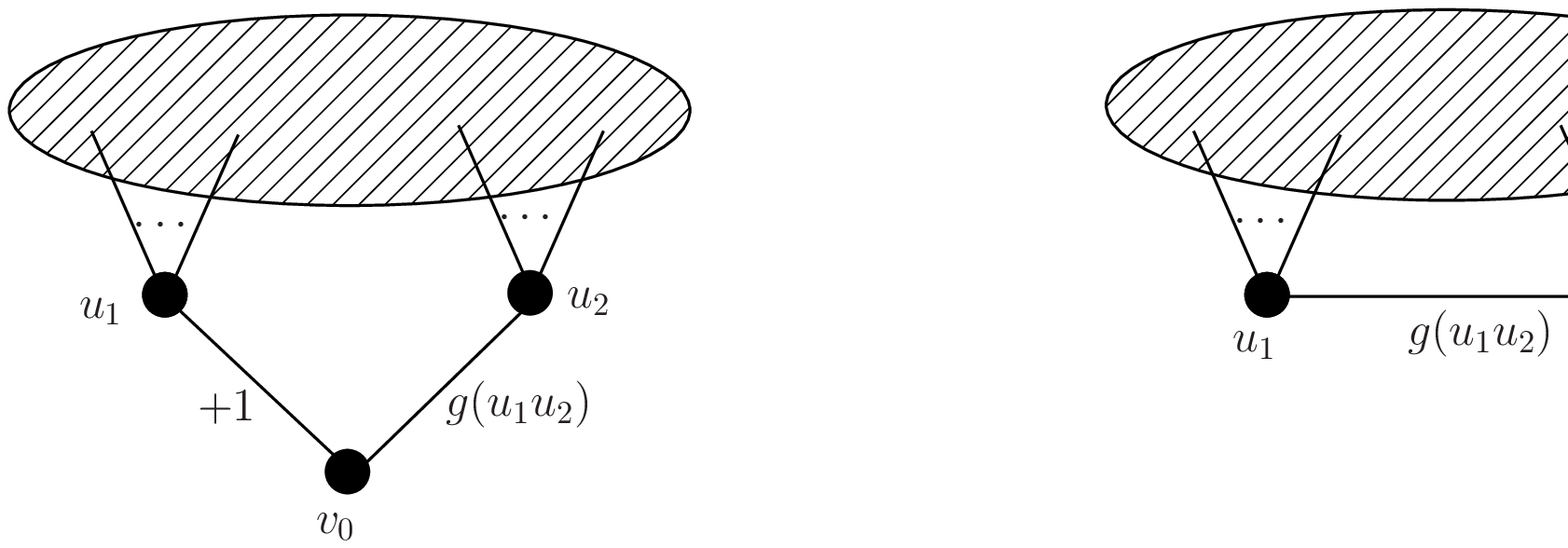}\\

\hspace{1 cm}(a) $G$ \hspace{4 cm} (b) $G'$ ($=G+u_1u_2-v_0$)

\caption{Graphs $G$ and $G'$ ($=G+u_1u_2-v_0$).}
\label{f5}
\end{figure}
\end{enumerate}
}

\proof (i). By the definition of $f$,
for any $v\in V(G)$, we have
$$
f(v)=
\ARRAY
{g(v), \quad &\mbox{if }v\in V(G)-\{v_0\};\\
2,  &\mbox{if }v=v_0.
}
$$
For any $uv\in E(G)-\{v_0u_1, v_0u_2, u_1u_2\}$ such that $f(uv)=1$,
we have $f(u)+f(v)=g(u)+g(v)\geq 2$.
For $v_0u_i$, $i\in \{1,2\}$, we have $f(v_0)+f(u_i)=2+g(u_i)\geq 2.$
Thus, $g\in \SEDF^0(G')$ implies $f\in \SEDF^0(G)$.

(ii). If $g(u_1u_2)=1$,
then by the definition of $f$, we have
$$
f(v)=
\ARRAY
{g(v), \quad &\mbox{if }v\in V(G)-\{v_0\};\\
2,  &\mbox{if }v=v_0.
}
$$
For any $uv\in E(G)-\{v_0u_1, v_0u_2\}$ such that $f(uv)=1$,
we have $f(u)+f(v)=g(u)+g(v)\geq 2$.
For $v_0u_i$, $i\in \{1,2\}$, we have $f(v_0)+f(u_i)=2+g(u_i)\geq 2.$
Thus, $f\in \SEDF^0(G)$ in this case.

If $g(u_1u_2)=-1$, then,
by the definition of $f$, we obtain
$$
f(v)=
\ARRAY
{g(v), \quad &\mbox{if }v\in V(G)-\{u_1,v_0\};\\
g(u_1)+2,  &\mbox{if }v=u_1;\\
0,  &\mbox{if }v=v_0.
}
$$
For any $uv\in E(G)-\{v_0u_1, v_0u_2\}$ such that $f(uv)=1$,
we have $f(u)+f(v)\geq g(u)+g(v)\geq 2$.
For the positive edge $v_0u_1$, we have $f(v_0)+f(u_1)=0+g(u_1)+2\geq 2.$
Thus, $f\in \SEDF^0(G)$ in this case.
\proofend

\section{$\gamma_s'(G)\le n+v_{odd}(G)/2$
\label{sec3}
}

For any graph $G$ and $f:E(G)\rightarrow \{1,-1\}$,
let $I_f(G)=\{v\in V(G): f(v)=0\}$.
We will prove the main result in this section
by applying the following result due to
Karami, Khodkar, and Sheikholeslami
\cite{karami2009improved}.

\begin{theo}[\cite{karami2009improved}]\label{th-kara}
For any simple graph $G$ of order $n$ with $v_{odd}(G)=0$,
there exists $f\in \SEDF^0(G)$ with
$I_f(G)\ne \emptyset$
and $f(u)\in \{0,2\}$ for all $u\in V(G)$.
\end{theo}

\prop
{pro3-1}
{For any simple graph $G$ of order $n$,
there is $f\in \SEDF^0(G)$ such that
$f(G)\le n+v_{odd}(G)/2$.}

\proof
It is sufficient to prove for the case when $G$ is connected.
It can be easily verified that the result holds
whenever $n\le 3$. Now assume that $n\ge 4$
and the result holds for any connected graph of order at most $n-1$.

If $G$ is not 2-connected,
by assumption, the result holds for each block of $G$.
Assume that $G$ has $k$ blocks, then by using Lemma~\ref{le2-2}
$k-1$ times, we will see the result holds for $G$.
If $G$ is $2$-connected and $\delta(G)=2$,
then the result also holds by assumption and
Lemma~\ref{le2-3}.

If $v_{odd}(G)=0$, then the result follows from
Theorem~\ref{th-kara}.
In the following, we assume that $G$ is $2$-connected
with $\delta(G)\ge 3$ and $v_{odd}(G)>0$.

For convenience, let
$k=v_{odd}(G)/2$ in the proof, where $k\ge 1$.
Let $U=\{u_1,u_2,\ldots,u_{2k}\}$ be the set of
odd vertices in $G$,
and let $G'$ be the graph obtained from $G$
by adding a new vertex $w$ and $2k$ new edges joining
$w$ to all vertices in $U$.
Clearly, $v_{odd}(G')=0$
and $E(G')=E(G)\cup \{wu_i:1\le i\le 2k\}$.

By Theorem~\ref{th-kara},
there exists
$g\in \SEDF^0(G')$ with $I_g(G')\ne \emptyset$
and $g(u)\in \{0,2\}$ for all $u\in V(G')$.
Thus $g(w)\in \{0,2\}$.
As $E(G')=E(G)\cup E_{G'}(w)$,
$g(G')=g(w)+g(E(G))$ holds.
Thus, we have the following conclusion.

\inclaim $g(E(G))=g(G')-g(w)$.
\def \cli {1}

Let $U_1$ be the set of vertices $u_i\in U$
with $g(wu_i)=+1$ and $U_2=U-U_1$.
As $g(w)=|U_1|-|U_2|$ and $|U_1|+|U_2|=d_{G'}(w)=2k$,
the following conclusion holds.

\inclaim $|U_1|=k+g(w)/2$.
\def \clii {2}

\begin{figure}[ht!]
\centering

\includegraphics[width=12cm]{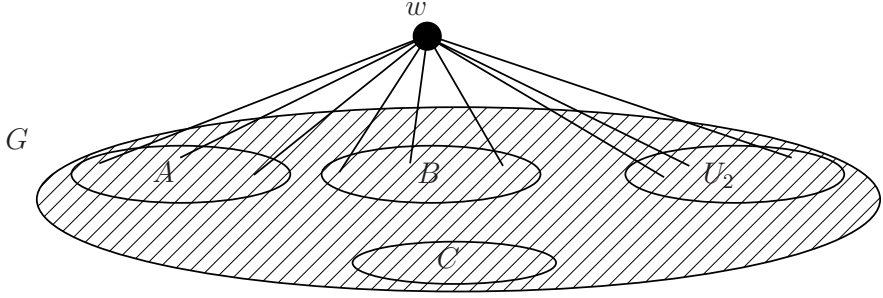}\\

\caption{$G=G'-w$, $N_{G'}(w)=A\cup B\cup U_2$ and
$B\cup C\subseteq I_g(G')$.}
\label{f22}
\end{figure}

$U_1$ is then partitioned into $A$ and $B$,
where $A$ is the set of $u_i\in U_1$ with $g(u_i)=2$.
Let $C$ be the set of vertices $v\in V(G)-U$
with $g(v)=0$ as shown in Figure~\ref{f22}.
Then $B\cup C\subseteq I_g(G')$.
Note that $w\in I_g(G')$ if and only if $g(w)=0$.
Thus the following claim holds.

\inclaim $|I_g(G')|\ge |B|+|C|+1-g(w)/2$.
\def \cliii {3}

By Theorem~\ref{th-kara}, we have
$g(G')=\frac 12 \sum_{u\in V(G')}g(u)=(n+1)-|I_{g}(G')|$.
Thus, the following conclusions follows from
Claims~\cli\ and~\cliii.

\inclaim $g(E(G))\le n-(|B|+|C|+g(w)/2)$.
\def \cliv {4}

Let $v$ be any vertex in $V(G)$.
As $\delta(G)\ge 3$, $d_{G'}(v)\ge 4$ holds.
Since $g(v)\in \{0,2\}$, $v$ is incident with some
edge $e\in E(G)$ with $g(e)=-1$.
Thus, there exists a subset $E_1$ of $E(G)$
with $g(e)=-1$ for all $e\in E_1$
such that each $v\in A\cup B$ is incident with some edge in $E_1$
(recall that $g(uw)=1$ for $u\in A\cup B=U_1$).
Let $E_1$ be a minimal one of such sets; note that $|E_1|\le |A|+|B|$.

Let $f:E(G)\rightarrow \{+1,-1\}$ be the function
defined by $f(e)=+1$ for all $e\in E_1$
and $f(e)=g(e)$ for all $e\in E(G)-E_1$.
It can be easily verified that
$f\in \SEDF^0(G)$ holds by the following facts:

\begin{enumerate}
\item For each $u_i\in A\cup B$,
we have $f(u_i)\ge g(u_i)-g(wu_i)+2=g(u_i)+1\ge 1$.

\item For each $u_i\in U_2$,
we have $f(u_i)\ge g(u_i)-g(wu_i)=g(u_i)+1\ge 1$.

\item For each $v\in V(G)-U$,
we have $f(v)\ge g(v)\ge 0$.

\item For each $e=v_1v_2\in E(G)$ with $f(e)=+1$,
if $e\in E_1$, then $f(v_1)+f(v_2)\ge
g(v_1)+g(v_2)+2\ge 2$;
if $e\in E(G)-E_1$, then $f(v_1)+f(v_2)\ge
g(v_1)+g(v_2)\ge 2$.
\end{enumerate}
By the definition of $f$ and Claim~\cliv, we have
\begin{eqnarray*}
\ f(G) & = & g(E(G))+2|E_1| \\
& \leq & n-(|B|+|C|+g(w)/2)+2(|A|+|B|) \\
& = & n+2|A|+|B|-|C|-g(w)/2.
\end{eqnarray*}
Thus, the following conclusion holds.

\inclaim $\gamma_s'(G)\le n+2|A|+|B|-|C|-g(w)/2$.
\def \clv {5}

Similarly, there exists a subset $E_2$ of $E(G)$
with $g(e)=-1$ for all $e\in E_2$
such that
each $v\in B\cup C$ is incident with some edge in $E_2$.
Let $E_2$ be a minimal one of such sets;
note that $|E_2|\le |B|+|C|$.

Let $f':E(G)\rightarrow \{+1,-1\}$  be the function
defined by $f'(e)=+1$ for all $e\in E_2$
and $f'(e)=g(e)$ for all $e\in E(G)-E_2$.
Again, it can be verified easily that
$f'\in \SEDF^0(G)$ holds by the following facts:
\begin{enumerate}
\item For each $u\in A$,
we have $f'(u)\ge g(u)-1=2-1=1$.

\item For each $u\in B$,
we have $f'(u)\ge g(u)-1+2=0-1+2=1$.

\item For each $u\in U_2$,
we have $f'(u)\ge g(u)+1\geq 0+1=1$.

\item For each $u\in C$,
we have $f'(u)\ge g(u)+2=0+2\geq 2$.

\item For each $u\in V(G)-U-C$,
we have $f'(u)\ge g(u)=2$.

\item For each $e=v_1v_2\in E(G)$ with $f'(e)=+1$,
we have $f(v_1)+f(v_2)\ge 1+1=2$.
\end{enumerate}
By the definition of $f'$ and Claim~\cliv, we have
\begin{eqnarray*}
\ f'(G) & = & g(E(G))+2|E_2| \\
& \leq & n-(|B|+|C|+g(w)/2)+2(|B|+|C|) \\
& = & n+|B|+|C|-g(w)/2.
\end{eqnarray*}
Thus, the following conclusion holds.

\inclaim $\gamma_s'(G)\le n+|B|+|C|-g(w)/2$.
\def \clvi {6}

By Claims~\clv\ and~\clvi,
$\gamma_s'(G)\le n+|A|+|B|-g(w)/2
=n+|U_1|-g(w)/2$ holds.
By Claim~\clii, we have $\gamma_s'(G)\le n+k
=n+v_{odd}(G)/2$.
\proofend

\section{
$\gamma_s'(G)\le n-2+v_{even}(G)$ when $v_{even}(G)>0$
\label{sec4}
}

In this section, the following exact values of $\gamma_s'(K_{m,n})$ will be used.

\begin{theo}[\cite{akbari2009signed}]\label{Kmn}
Let $m$ and $n$ be two positive integers, where $m\leq n$. Then:
\begin{enumerate}
\item If $m$ and $n$ are even, then $\gamma_s'(K_{m,n})=\min\{2m, n\}$.

\item If $m$ and $n$ are odd, then $\gamma_s'(K_{m,n})=\min\{2m-1, n\}$.

\item If $m$ is even and $n$ is odd, then $\gamma_s'(K_{m,n})=\min\{3m, \max\{2m, n+1\}\}$.

\item If $m$ is odd and $n$ is even, then $\gamma_s'(K_{m,n})=\min\{3m-1, \max\{2m, n\}\}$.
\end{enumerate}
\end{theo}

In the proof of the part $(b)$ of Theorem \ref{th2},
we shall need the parts (i) and (iii) of Theorem \ref{Kmn}.
In these two cases, actually Akbari et al. proved that there exists
$f\in \SEDF^0(G)$ such that $f(K_{m,n})=\gamma_s'(K_{m,n}).$

\prop
{pro41}
{For any simple graph $G$ of order $n$, if $v_{even}(G)>0$,
then there is an $f\in \SEDF^0(G)$ such that $f(G)\le n-2+v_{even}(G)$,
and thus $\gamma_s'(G)\le n-2+v_{even}(G)$.}

\proof
When $v_{even}(G)=0$,
by Theorem 7 in \cite{karami2009improved},
there exists
$f\in \SEDF^0(G)$ with $f(G)\le |V(G)|-1$.
So it is sufficient to prove the case when $G$ is connected.
It can be easily verified that the result holds
whenever $n\le 3$. Now assume that $n\ge 4$
and the result holds for any graph of order at most $n-1$.
By Lemma \ref{le2-2},
we only need to prove the result for
$2$-connected graphs.
Let $v_{even}(G)=t\geq 1$, and let $W=\{w_1, w_2, \ldots w_t\}$ be the set of all even vertices.

\inclaim If $W$ is an independent set,
then there exists $f\in \SEDF^0(G)$ such that $f(G)\le n-2+v_{even}(G)$.
\def \cli {1}

Assume that $w_1$ has the minimum degree among all elements in $W$.
Let $d_G(w_1)=2s$, $s\geq 1$, and assume
that $N_G(w_1)=\{u_1, u_2, \ldots, u_{2s}\}$.
Consider $G'=G-w_1$. Since $G$ has no cut vertex, $G'$ is connected.
Clearly, $|V(G')|=n-1$ and $v_{odd}(G')=n-t-2s$.

{\bf Case 1.1.} $n-t-2s\geq 2$, i.e., $G'$ is not an Eulerian graph.

In this case, $G'$ can be decomposed into $(n-t-2s)/2$ trails $T_1$, \ldots, $T_{(n-t-2s)/2}$,
and the endpoints of these $(n-t-2s)/2$ trails correspond to all odd vertices of $G'$.
Now, we define the function $f_1: E(G)\rightarrow \{1,-1\}$ as follows:
\begin{enumerate}
\item for each $T_i$, $1\leq i\leq (n-t-2s)/2$, starting with $+1$,
we assign $+1$ and $-1$ to the edges of $T_i$ alternatively.
When the trail has even number of edges, we change the value of the last edge to $+1$;

\item for each edge $w_1u_i$, $1\leq i\leq 2s$, we set $f_1(w_1u_i)=+1$; and

\item for any $w_i$, $2\leq w_i\leq t$, if the weight of $w_i$ till now is $0$, then we choose
      any negative edge incident to $w_i$ and change it to a positive one.
\end{enumerate}
It follows from the construction that
$$f_1(G)\leq 2\cdot \frac{n-t-2s}{2}+2s+2(t-1)=n-2+t=n-2+v_{even}(G).$$
Next, after Step (i), the weight of any vertex in $V(G)-W-N_G(w_1)$ is at least 1,
and the weight of vertices in $(W-\{w_1\})\cup N_G(w_1)$ is 0 or at least 2.
After Step (ii), $f(w_1)$ is $2s$; the weight of vertices in $N_G(w_1)$ has increased by 1,
and others remain unchanged.
Finally, after Step (iii), all the vertices in $W-\{w_1\}$
are of the weight at least 2.

Hence $f_1\in \SEDF^0(G)$, and $\gamma_s'(G)\leq f_1(G)\leq n-2+v_{even}(G)$.

{\bf Case 1.2.} $n-t-2s=0$, i.e., $G'$ is an Eulerian graph, and $2s\geq 4$.

Because $G'$ is an Eulerian graph, so it has an Eulerian circuit.
Now we define the function $f_2: E(G)\rightarrow \{1,-1\}$ as follows:
\begin{enumerate}
\item for a fixed Eulerian circuit of $G'$, starting from the vertex $u_1$, walking along the Eulerian circuit,
we assign $+1$ and $-1$ alternatively starting with $+1$;

\item we set $f_2(w_1u_i)=1$, $1\leq i\leq 2s$ if $|E(G')|$ is even;
otherwise, if $|E(G')|$ is odd, we set $f_2(w_1u_1)=-1$ and $f_2(w_1u_i)=1$, $2\leq i\leq 2s$; and

\item for any $w_i$, $2\leq w_i\leq t$, choose any negative edge incident to $w_i$ and change it to a positive one.
\end{enumerate}

After Step (i), if $G'$ has an even number of edges, then each vertex in $G'$ has weight 0;
if $G'$ has an odd number of edges, then $f_2(u_1)=2$, and all other vertices have weight 0.
Next, after Step (ii), all vertices in $N_G(w_1)$ have weight 1, and $f_2(w_1)\geq 2s-2\geq 2$, and
all vertices in $W-\{w_1\}$ have weight 0.
Finally, after Step (iii), all vertices in $W-\{w_1\}$ have weight 2, and others do not decrease.

Hence $f_2\in \SEDF^0(G)$. Recall that in this case $n-t-2s=0$, we have
$$\gamma_s'(G)\leq f_2(G)=0+2s+2(t-1)=2s+2t-2=n+t-2=n-2+v_{even}(G)$$
when $|E(G')|$ is even, and
$$\gamma_s'(G)\leq f_2(G)=1+(2s-2)+2(t-1)=2s+2t-3=n+t-3=n-3+v_{even}(G)$$
when $|E(G')|$ is odd.

{\bf Case 1.3.} $n-t-2s=0$ and $2s=2$.

In this case, if $n$ is odd, then $G=K_{2, n-2}$.
Then, if $n\geq 5$, there exists $f_3\in \SEDF^0(G)$ such that
$f_3(G)=\gamma_s'(G)=\gamma_s'(K_{2, n-2})=\min\{6, \max\{4, n-1\}\}$.
Hence there exists $f_3\in \SEDF^0(G)$ such that
\begin{displaymath}
\gamma_s'(G) = f_3(G) = \left\{ \begin{array}{ll}
2, & \textrm{if $n=3$,} \\
4, & \textrm{if $n=5$,} \\
6, & \textrm{if $n\geq 7$.}
\end{array} \right.
\end{displaymath}
Therefore $\gamma_s'(G)=f_3(G)\leq n-2+v_{even}(G)$ holds.

If $n$ is even, then $G=K_{2, n-2}+u_1u_2$.
There exists $f\in \SEDF^0(K_{2, n-2})$ such that
$f(K_{2, n-2})=\gamma_s'(K_{2, n-2})=\min\{4, n-2\}$, $n\geq 4$.
Now we extend $f$ by assigning $+1$ to $u_1u_2$, thus obtain $f_4\in \SEDF^0(G)$ such that
$f_4(G)=\gamma_s'(K_{2, n-2})+1=\min\{4, n-2\}+1$.
That is, there exists $f_4\in \SEDF^0(G)$ such that
\begin{displaymath}
\gamma_s'(G)\leq f_4(G) = \left\{ \begin{array}{ll}
3, & \textrm{if $n=4$,} \\
5, & \textrm{if $n\geq 6$.}
\end{array} \right.
\end{displaymath}
Therefore $\gamma_s'(G)\leq f_4(G)\leq n-2+v_{even}(G)$ holds.

\inclaim
If $W$ is not an independent set,
then there is an $f\in \SEDF^0(G)$ such that $f(G)\le n-2+v_{even}(G)$.
\def \clii {2}

In this case, we find a maximal matching $M$ in $G[W]$. Assume that $|M|=p$ and
consider the graph $G''=G-M$, with $n$ vertices
and $v_{even}(G'')=t-2p$.

{\bf Case 2.1.} $t-2p\geq 1$.

Since $M$ is maximal in $G[W]$, the $t-2p$ even vertices in $G''$ form an independent set.
By Claim~\cli,
there is an $f_5\in \SEDF^0(G'')$ such that $f_5(G'')\le n-2+v_{even}(G'')$.
We now extend $f_5$ by adding $M$ to $G''$ and letting each edge in $M$ be a positive edge.
Thus we obtain $f_5'\in\SEDF^0(G)$ such that $\gamma_s'(G)\leq f_5'(G)= f_5(G'')+p\leq n-2+t-2p+p=n-2+t-p<n-2+v_{even}(G)$.

{\bf Case 2.2.} $t-2p=0$, i.e., $M$ is a perfect matching of $G[W]$.

In this subcase, $v_{odd}(G'')=n$.
Karami et al. \cite{karami2009improved} proved that for
a graph $G$ with $n$ vertices in which each vertex is of odd degree, there
exists $f\in \SEDF^0(G)$ such that $\gamma_s'(G)\leq n-1$. So
there is $f_6\in \SEDF^0(G'')$ such that $f_6(G'')\le n-1$.
We now extend $f_6$ by adding $M$ to $G''$ and letting each edge in $M$ be a positive edge.
Thus we obtain $f_6'\in\SEDF^0(G)$ such that
$\gamma_s'(G)\leq f_6'(G)= f_6(G'')+p\leq n-1+p=n-1+\frac{v_{even}(G)}{2}\leq n-2+v_{even}(G)$.

Thus, Claim~\clii\ holds and the proof is complete.
\proofend

\begin{cor}
Conjecture $\ref{conj1}$ holds for the case $v_{even}(G)=1$.
\end{cor}

Now we prove Theorem \ref{th2}.

\noindent
{\bf Proof of Theorem \ref{th2}.}
From Propositions \ref{pro3-1} and \ref{pro41}, we can see that
$\gamma_s'(G)\le n+v_{odd}(G)/2$, and $\gamma_s'(G)\le n-2+v_{even}(G)$ when $v_{even}(G)>0$.

So when $v_{even}(G)>0$, we have
$$3\gamma_s'(G)\leq 2(n+v_{odd}(G)/2)+(n-2+v_{even}(G))=3n+v_{odd}(G)+v_{even}(G)-2=4n-2,$$
and hence $\gamma_s'(G)\le (4n-2)/3$.

When $v_{even}(G)=0$, i.e., $v_{odd}(G)=n$, it was proved in \cite{karami2009improved} that
$\gamma_s'(G)\leq n-1$. Hence $\gamma_s'(G)\le (4n-2)/3$ also holds.
{\hfill$\Box$}

\section{Conjecture~\ref{conj1} for $v_{even}(G)=2$
\label{sec5}
}

\prop
{pro5-1}
{For any simple graph $G$ of order $n$,
if $v_{even}(G)=2$, then
there is an $f\in \SEDF^0(G)$ such that $f(G)\le n-1$,
and thus $\gamma_s'(G)\le n-1$.
}

\proof
It can be easily verified that the result holds
whenever $n\le 3$. So assume that $n\ge 4$
and the result holds for any graph of order at most $n-1$.

Let $G$ be a simple graph of order $n$ with
$v_{even}(G)=2$, and let $w_1, w_2$ be the two vertices of even degree.

\inclaim Proposition~\ref{pro5-1} holds for $G$
when it is disconnected.
\def \clmi {1}

Assume that $G_1,G_2,\ldots,G_k$
are the components of $G$,
where $k\ge 2$.
Then $v_{even}(G_i)\le 2$ and $|V(G_i)|\le n-1$
for all $i=1,2,\ldots,k$.
For any $G_i$, where $1\le i\le k$,
if $v_{even}(G_i)=0$,
by the proof in \cite{karami2009improved},
there exists
$f_i\in \SEDF^0(G_i)$ with $f_i(G_i)\le |V(G_i)|-1$;
otherwise,
by the assumption above and
Proposition~\ref{pro41},
there exists
$f_i\in \SEDF^0(G_i)$ such that $f_i(G_i)\le |V(G_i)|-1$.
\\
Let $f$ be the mapping $E(G)\rightarrow \{1,-1\}$ defined
by $f|_{E(G_i)}=f_i$ for all $1\le i\le k$.
It is easy to see that
$f\in \SEDF^0(G)$ with
$$
f(G)=\sum_{i=1}^k f(G_i)\le \sum_{i=1}^k (|V(G_i)|-1)
= n-k\le n-2.
$$
Thus, Claim~\clmi\ holds.

\inclaim Proposition~\ref{pro5-1} holds for $G$
when $d(w_i)=2$ for some
$i\in \{1,2\}$.
\def \clmii {2}

Assume that $d(w_1)=2$.
Let $N(w_1)=\{u_1, u_2\}$.

If $u_1u_2\in E(G)$, then consider the graph
$G-u_1u_2-w_1$. $G-u_1u_2-w_1$ is a simple graph of order $n-1$ and $v_{even}(G-u_1u_2-w_1)=1$.
By Proposition \ref{pro41},
there exists $g\in \SEDF^0(G-u_1u_2-w_1)$ such that $g(G-u_1u_2-w_1)\le (n-1)-2+1=n-2$.
Then, by Lemma \ref{le2-3} (i),
there exists $f\in \SEDF^0(G)$ such that
$f(G)=g(G-u_1u_2-w_1)+1\leq n-1$.
If $u_1u_2\notin E(G)$,
then similarly, by applying Proposition \ref{pro41}
and Lemma \ref{le2-3} (ii),
we can show that there exists $f\in \SEDF^0(G)$ with
$f(G)\leq n-1$.

Thus, Claim~\clmii\ holds.

According to Claims 1 and 2, in the following,
we may assume that $G$ is connected and $d(w_i)\geq 4$ for $i=1, 2$.

Let $N_0=N_G(w_1)\cap N_G(w_2)$,
$N_1=N_G(w_1)-N_G(w_2)-\{w_2\}$,
$N_2=N_G(w_2)-N_G(w_1)-\{w_1\}$, and $N_3=V(G)-(N_0\cup N_1\cup N_2\cup \{w_1, w_2\})$.
Set $n_i=|N_i|$ for $0\leq i\leq 3$. Then $n_0+n_1+n_2+n_3=v_{odd}(G)=n-2$.

{\bf Case 1.} $w_1w_2\notin E(G)$.

We divide this case into two subcases, depending on
whether $w_1$ and $w_2$ share a neighbour or not,
i.e., $n_0\geq 1$ or $n_0=0$.

{\bf Case 1.1.} $w_1w_2\notin E(G)$ and $n_0\geq 1$.

Consider the graph $G'=G-w_1$.
Note that $v_{odd}(G')=n_2+n_3$, and so $G'$ can be decomposed into $(n_2+n_3)/2$ trails
$\{T_1, T_2, \ldots, T_{(n_2+n_3)/2}\}$, whose endpoints corresponds to
all odd vertices of $G'$.

If $T_i$ has odd length, we assign $+1$ and $-1$ alternatively to the edges of $T_i$,
starting and ending with $+1$; this weight assignment in an odd trail
is called a {\it proper assignment}.
When $T_i$ has even length $t$, there are exactly $\frac{t}{2}$ vertices on $T_i$,
each of which can naturally divide $T_i$ into two subtrails with odd length.
We call these $\frac{t}{2}$ vertices \emph{good}.
For each $T_i$ with even length, choose a good vertex $u_i$
of $T_i$.
We can assign $+1$ and $-1$ alternatively
to edges in the two subtrails of $T_i$ divided by $u_i$
such that both starting and ending edges in each subtrail
are assigned $+1$.
This weight assignment of edges in an even trail
$T_i$ is called
a \emph{proper assignment with respect to $u_i$}.
Let $\mathbb{T}_1=\{T_1, T_2, \ldots, T_{(n_2+n_3)/2}\}$.

\inclaim In Case 1.1,
Proposition \ref{pro5-1} holds for $G$ when
there is at least one trail of odd length in $\mathbb{T}_1$.
\def \clmiii {3}

Assume that there is at least one trail of odd length in $\mathbb{T}_1$.
For each $T_i\in \mathbb{T}_1$ with even length,
let $u_i$ be a good vertex of $T_i$.
We define a function $f_1: E(G)\rightarrow \{1,-1\}$ as follows:
each odd trial $T_i\in \T_1$ is equipped with
a proper assignment,
and each even trail $T_i\in \T_1$ is equipped with
a proper assignment with respect to $u_i$.
Then we assign $+1$ to each edge incident to $w_1$.
If the weight of $w_2$ till now is $0$, we choose any negative edge incident to $w_2$ and change it to a positive one.

Now we have $f_1(w_1)=d(w_1)\geq 4$, $f_1(w_2)\geq 2$, and $f_1(u)\geq 1$ for each $u\in V(G)-\{w_1, w_2\}$. So $f_1\in \SEDF^0(G)$
and hence $\gamma_s'(G)\leq f_1(G)\leq 1+2(\frac{n_2+n_3}{2}-1)+n_1+n_0+2=n-1$.
Thus Claim~\clmiii\ holds.

\inclaim In Case 1.1, if all trails in $\T_1$ have even length, then either $w_2$
or some vertex $x\in N_0$ is a good vertex of some
trail $T_j\in \T_1$.

\def \clmiv {4}

Assume that all trails in $\mathbb{T}_1$ have even length.
Then, some edge $w_2x$, where $x\in N_0$,
must be in some $T_j\in \T_1$.
Obviously, either $w_2$ or $x$ is a good vertex in $T_j$.
Thus Claim~\clmiv\ holds.

\inclaim In Case 1.1,  Proposition \ref{pro5-1} holds for $G$
when all trails in $\mathbb{T}_1$ have even length.
\def \clmv {5}

Assume that all trails in $\mathbb{T}_1$ have even length.
By Claim~\cliv, either $w_2$ or some vertex $x\in N_0$
is a good vertex of some trail $T_j\in \T_1$.

If $w_2$ is good, we define a function $f_2: E(G)\rightarrow \{1,-1\}$ as follows.
We equip $T_j$ with the proper assignment with respect to $w_2$,
and for each $T_i\in \mathbb{T}_1-\{T_j\}$, we equip $T_i$ with a proper assignment (with respect to any good vertex).
Then we assign $+1$ to each edge incident to $w_1$.

Now we have $f_2(w_1)=d(w_1)\geq 4$, $f_2(w_2)\geq 2$, and $f_2(u)\geq 1$ for each $u\in V(G)-\{w_1, w_2\}$.
So $f_2\in \SEDF^0(G)$ and hence $\gamma_s'(G)\leq f_2(G)= 2\cdot\frac{n_2+n_3}{2}+n_1+n_0=n-2$.

If $x$ is good, we define a function $f_3: E(G)\rightarrow \{1,-1\}$ as follows.
We equip $T_j$ with the proper assignment with respect to $x$,
and for each $T_i\in \mathbb{T}_1-\{T_j\}$,
we equip $T_i$ with a proper assignment (with respect to any good vertex).
Then we assign $-1$ to $w_1x$ and $+1$ to any other edge incident to $w_1$.
If the weight of $w_2$ till now is $0$, we choose any negative edge incident
to $w_2$ and change it to a positive one.

Now we have $f_3(w_1)=d(w_1)-2\geq 2$, $f_3(w_2)\geq 2$, $f_3(x)\geq 2-1=1$,
and $f_3(u)\geq 1$ for each $u\in V(G)-\{w_1, w_2, x\}$.
So $f_3\in \SEDF^0(G)$ and hence $\gamma_s'(G)\leq f_3(G)\leq 2\cdot\frac{n_2+n_3}{2}+(n_1+n_0-2)+2=n-2$.
Thus Claim~\clmv\ holds.

By Claims~\clmiii\ and~\clmv,
Proposition~\ref{pro5-1} holds for $G$ in Case 1.1.

{\bf Case 1.2.} $w_1w_2\notin E(G)$ and $n_0= 0$.

\inclaim Proposition~\ref{pro5-1} holds for $G$ in Case 1.2.
\def \clmvi {6}

Choose edges $e_1, e_2$ incident to $w_1$
and edges $e_3,e_4$ incident to $w_2$.
\\
\begin{figure}[ht!]
\centering
\includegraphics[width=12cm]{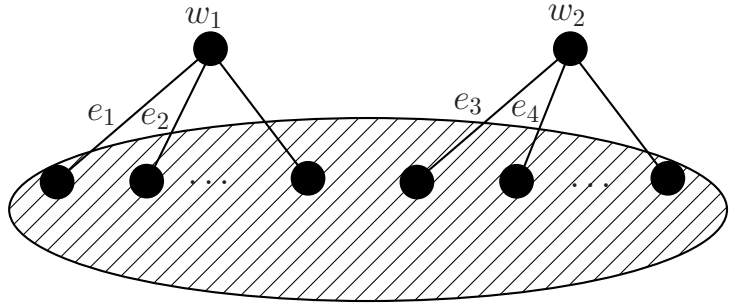}\\
\caption{Case 1.2.}
\label{f23}
\end{figure}
\\ As $n_0=0$, $N_0=N(w_1)\cap N(w_2)=\emptyset$.
Thus, by the condition that $d(w_i)\ge 4$ for both $i=1,2$,
we have $n\ge 2+4\cdot 2=10$.
Let $G''$ denote the graph $G-\{e_1,e_2,e_3,e_4\}$.
Observe that $v_{odd}(G'')=v_{odd}(G)-4=n-6>0$.
Thus $E(G'')$ can be decomposed into $t=(n-6)/2$ trails,
say $T_1,T_2,\ldots, T_t$.
We now define the function $f_4:E(G)\rightarrow \{1,-1\}$ as follows:
\begin{itemize}
\item $f_4(e_i)=1$ for $i=1,2,3,4$; and
\item each odd trial $T_i$ is equipped with
a proper assignment,
and each even trail $T_i$ is equipped with
a proper assignment
with respect to some good vertex.
\end{itemize}
Observe that $f_4(w_i)\ge 2$ for $i=1,2$
and $f_4(u)\ge 1$ for all $u\in V(G)-\{w_1,w_2\}$.
Thus $f_4\in \SEDF^0(G)$.
Also note that
$
f_4(G)\le 2t+4=2(n-6)/2+4=n-2.
$
So Claim~\clmvi\ holds.

{\bf Case 2.} $w_1w_2\in E(G)$.

Similarly as in Case 1, we divide this case into two subcases,
depending on whether $w_1$ and $w_2$ share a neighbour or not.

{\bf Case 2.1.} $w_1w_2\in E(G)$ and $n_0\geq 1$.

\inclaim Proposition~\ref{pro5-1} holds for $G$ in Case 2.1.\def \clmvii {7}

Consider the graph $G'=G-w_1$.
Note that $v_{odd}(G')=n_2+n_3+1$,
and so $G'$ can be decomposed into $(n_2+n_3+1)/2$ trails
$\{T_1, T_2, \ldots, T_{(n_2+n_3+1)/2}\}$ whose endpoints
correspond to all odd vertices of $G'$.
Let $\mathbb{T}_2=\{T_1, T_2, \ldots, T_{(n_2+n_3+1)/2}\}$.

Case 2.1 is now divided into two subcases.

{\bf Case 2.1.1.} Some trail in $\mathbb{T}_2$ has an odd length.

We define a function $g_1: E(G)\rightarrow \{1,-1\}$ as follows.
Each trail $T_i\in \mathbb{T}_2$ of odd length
is equipped with a proper assignment,
and each trail $T_i\in \mathbb{T}_2$
of even length
is equipped with a proper assignment with respect to
some good vertex of $T_i$.
Then we assign $+1$ to each edge incident to $w_1$.

Now we have $g_1(w_1)=d(w_1)\geq 4$, $g_1(w_2)\geq 2$, and $g_1(u)\geq 1$ for each $u\in V(G)-\{w_1, w_2\}$. So $g_1\in \SEDF^0(G)$
and hence $\gamma_s'(G)\leq g_1(G)\leq 1+2(\frac{n_2+n_3+1}{2}-1)+n_1+n_0+1=n-1$.

{\bf Case 2.1.2.}
All trails in $\mathbb{T}_2$ have even length.

Choose any $x\in N_0$ and assume that
$w_2x$ is an edge in $T_1$.
Then, either $w_2$ or $x$ is good in $T_1$.
Let $u_1=w_2$ if $w_2$ is good in $T_1$,
and $u_1=x$ otherwise.
For any $i=2,3,\ldots,(n_2+n_3+1)/2$,
let $u_i$ be any good vertex of $T_i$.

We define a function $g_2: E(G)\rightarrow \{1,-1\}$ as follows.
We first equip each $T_i$ with a proper assignment
with respect to $u_i$.
Then, we assign $-1$ to $w_1u_1$, and finally, we assign $+1$ to any
other edge incident with $w_1$.

If $u_1=w_1$, then $g_2(w_1)=d(w_1)-2\geq 2$,
$g_2(w_2)\geq 2-1=1$,
and $g_2(u)\geq 1$ for each $u\in V(G)-\{w_1, w_2\}$.
So $g_2\in \SEDF^0(G)$ and hence $\gamma_s'(G)\leq g_2(G)\leq 2\cdot\frac{n_2+n_3+1}{2}+n_1+n_0-1=n-2$.

If $u_1=x$, then
$g_2(w_1)=d(w_1)-2\geq 2$,
$g_2(w_2)\geq 1$, $g_2(x)\geq 2-1=1$,
and $g_2(u)\geq 1$ for each $u\in V(G)-\{w_1, w_2, x\}$.
So $g_2\in \SEDF^0(G)$ and hence $\gamma_s'(G)\leq g_2(G)\leq 2\cdot\frac{n_2+n_3+1}{2}+(n_1+n_0-2+1)=n-2$.

Hence Claim~\clmvii\ holds.

{\bf Case 2.2.} $w_1w_2\in E(G)$ and $n_0=0$.

\inclaim Proposition~\ref{pro5-1} holds for $G$ in Case 2.2.
\def \clmviii {8}

Choose $x_1\in N_1$ and $x_2\in N_2$
and consider the graph $G'''=G-\{e_0,e_1,e_2\}$, where
$e_0=w_1w_2$, $e_1=w_1x_1$, and $e_2=w_2x_2$.
\\
\begin{figure}[ht!]
\centering
\includegraphics[width=9.5cm]{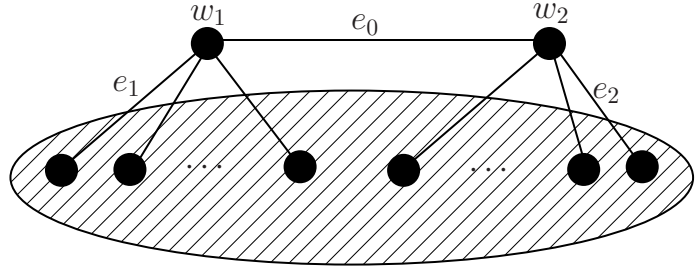}\\
\caption{Case 2.2.}
\label{f24}
\end{figure}
\\
As $n_0=0$, $N_0=N(w_1)\cap N(w_2)=\emptyset$.
Thus, by the condition that
$d(w_i)\ge 4$ for both $i=1,2$,
we have $n\ge 2+3\cdot 2=8$.
Observe that $v_{odd}(G''')=v_{odd}(G)-2=n-4>0$.
Thus $E(G''')$ can be decomposed into $t=(n-4)/2$ trails,
say $T_1,T_2,\ldots, T_t$.
We now define a function $g_3:E(G)\rightarrow \{1,-1\}$ as follows:
\begin{itemize}
\item $g_3(e_i)=1$ for $i=0,1,2$; and
\item each odd trial $T_i$ is equipped with
a proper assignment
and each even trail $T_i$ is equipped with
a proper assignment
with respect to some good vertex.
\end{itemize}
Observe that $g_3(w_i)\ge 2$ for $i=1,2$,
and $g_3(u)\ge 1$ for all $u\in V(G)-\{w_1,w_2\}$.
Thus $g_3\in \SEDF^0(G)$.
Also note that
$g_3(G)\le 2t+3=2(n-4)/2+3=n-1$,
and so Claim~\clmviii\ holds,
which eventually finishes the proof.
\proofend

\vspace{0.3 cm}

Note that Theorem~\ref{th3} follows directly from
Proposition~\ref{pro41} for the case $v_{even}(G)=1$
and from Proposition~\ref{pro5-1} for the case $v_{even}(G)=2$.

\section{Concluding remarks
\label{sec6}
}

Karami et al. \cite{karami2009improved}
proved Conjecture~\ref{conj1} for the two cases
$v_{odd}(G)=0$ or $n$
by showing the existence of $f\in \SEDF^0(G)$
with $f(G)\le n-1$.
In the proof of Propositions~\ref{pro3-1}, \ref{pro41}
and \ref{pro5-1},
all defined members in $\SEDF(G)$
also belong to $\SEDF^0(G)$.
Therefore, we believe Conjecture~\ref{conj1} can be strengthened
to the following one.

\begin{con}\label{conj2}
For any simple graph $G$ of order $n$,
there exists  $f\in \SEDF^0(G)$
with $f(G)\le n-1$.
\end{con}

In 2005, Xu \cite{xu2005edge} proved the following sharp lower bound of $\gamma_s'(G)$.

\begin{theo}[\cite{xu2005edge}]\label{low1}
Let $G$ be a graph with $n$ vertices, $m$ edges and $\delta(G)\geq 1$. Then
$\gamma_s'(G)\geq n-m$.
\end{theo}

Then Karami et al. \cite{karami2008some} characterized all simple connected graphs $G$ for which $\gamma_s'(G)= n-m$.
These graphs all have many vertices of degree 1. If we restrict graphs to have higher connectivity
or larger minimum degree, a better lower bound can be expected. So we raise the following conjecture.

\begin{con}\label{conj3}
Let $G$ be a $2$-connected graph with $n$ vertices and $m$ edges, and
without two adjacent degree $2$ vertices. Then
$\gamma_s'(G)\geq 2n-m$.
\end{con}

If the conjecture above is correct, then the lower bound is also sharp.
For example, $\gamma_s'(K_4-e)=3=2n-m$.

Now we show more examples that the bound in Conjecture \ref{conj3} is reachable.
Let $G$ be a $2$-connected Hamiltonian graph
with $\delta(G)\ge 3$, $V(G)=\{v_1,v_2,\ldots,v_n\}$ and size $m$.
Suppose $C$ is one of its Hamiltonian cycles.

The \emph{triangulation} of a graph $H$, denoted by $T(H)$,
is the graph obtained from $H$ by changing each edge $uv$ of $H$
into a triangle $uwv$, where $w$ is a new vertex associated with $uv$.
Let $G'=T(G-E(C))+E(C)$, that is, the graph obtained from $T(G-E(C))$ by adding all the edges in the Hamiltonian cycle $C$.
Then the order of $G'$ is $m$ and
the size of $G'$ is $3m-2n$.

Observe that $G'$ is 2-connected and
does not have two adjacent degree 2 vertices.
Consider a function $f:E(G')\rightarrow \{1,-1\}$, where
$f(e)=1$ if $e\in E(G)$, and $f(e)=-1$ otherwise.
Then
$$
f(G')=|E(G)|-|E(G')-E(G)|=m-2(m-n)=2n-m=2|V(G')|-|E(G')|.
$$

\begin{figure}[ht!]
\centering

\includegraphics[width=7cm]{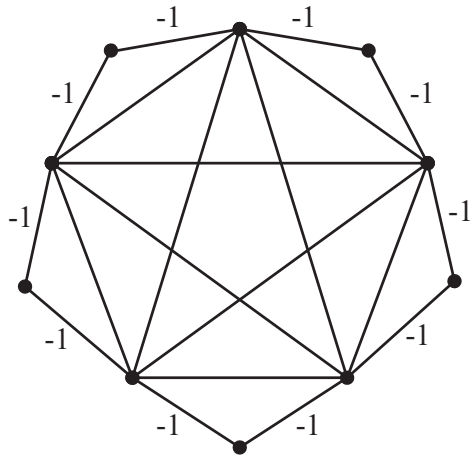}\\

\caption{$G'=T(G-E(C))+E(C)$, where $G=K_5$.}
\label{example-conj3}
\end{figure}

By the definition of $f$,
$f_{G'}(v_i)=2$ for each $i=1,2,\ldots,n$
whereas $f(u)=-2$ for each $u\in V(G')-V(G)$.
Thus, for each $e=uv\in E(G)$, we have $f(e)=1$
and $f(u)=f(v)=2$, whereas
for each $e=uv\in E(G')-E(G)$, we have $f(e)=-1$
and $f(u)+f(v)=0$.
Thus $f\in \SEDF(G)$. The graph shown in Figure \ref{example-conj3}
is an example of $G'$ when $G=K_5$ (edges without a sign in the figure receive sign $+1$).

\section*{Acknowledgements}

This paper is supported by NSFC (Nos. 11701401 and 11971346).
The authors would like to thank two anonymous referees for their very
helpful comments.
This research was conducted during the
second and third authors' visit at National Institute of Education, Nanyang Technological University, Singapore.
They are grateful to staff in NIE/NTU for their help during their visits.

\newcommand \AM[1]
{{\it Ann. of Math.} {\bf #1}}

\newcommand \AJC[1]
{{\it Australasian J. Combin.} {\bf #1}}

\newcommand \ARS[1]
{{\it ARS Combinatoria} {\bf #1}}

\newcommand \CJM[1]
{{\it Canadian J. Math.} {\bf #1}}

\newcommand \COMBIN[1]
{{\it Combinatorica} {\bf #1}}

\newcommand \CPC[1]
{{\it Combin. Probab. Comput.} {\bf #1}}

\newcommand \DAM[1]
{{\it Discrete Appl. Math.} {\bf #1}}

\newcommand \DM[1]
{{\it Discrete Math.} {\bf #1}}

\newcommand \EJC[1]
{{\it Electronic J. of Combin. Theory} {\bf #1}}

\newcommand \EUJC[1]
{{\it European J. of Combin. Theory} {\bf #1}}

\newcommand \GC[1]
{{\it Graphs Combin.} {\bf #1} }

\newcommand \JAMS[1]
{{\it J. Amer. Math. Soc.} {\bf #1}}

\newcommand \JCT[1]
{{\it J. Combin. Theory} {\bf #1}}

\newcommand \JCTA[1]
{{\it J. Combin. Theory Ser. A} {\bf #1}}

\newcommand \JCTB[1]
{{\it J. Combin. Theory Ser. B} {\bf #1}}

\newcommand \JCMCC[1]
{{\it J. Combin. Math. Combin. Comput.} {\bf #1}}

\newcommand \JGT[1]
{{\it J. Graph Theory} {\bf #1} }

\newcommand \JG[1]
{{\it J. Geom.} {\bf #1} }

\newcommand \LAA[1]
{{\it Linear Algebra Appl.} {\bf #1} }

\newcommand \LMA[1]
{{\it Linear Multilinear Algebra} {\bf #1} }

\newcommand \PCPS[1]
{{\it Proc. Cambridge Philos. Soc.} {\bf #1} }

\newcommand \TAMS[1]
{{\it Trans. Amer. Math. Soc.} {\bf #1}}

\newcommand \SIAMD[1]
{{\it SIAM J. Discrete Math.} {\bf #1}}

\newcommand \UM[1]
{{\it Util. Math.} {\bf #1}}

\begin{thebibliography}{10}

\bibitem{akbari2009signed}
S. Akbari, S. Bolouki, P. Hatami, and M. Siami,
On the signed edge domination number of graphs,
\DM{309} 
(3) (2009), 587--594.

\bibitem{akbari2014some}
S. Akbari, H. Esfandiari, E. Barzegary, and S. Seddighin,
Some bounds for the signed edge domination number of a graph,
\AJC{58} 
(2014), 60--66.

\bibitem{chen2008lower}
W. Chen and E. Song,
Lower bounds on several versions of signed domination number,
\DM{308} 
(10) (2008), 1837--1846.

\bibitem{dunbar1995signed}
J. Dunbar, S. Hedetniemi, M.A. Henning, and P.J. Slater,
Signed domination in graphs,
in: Y. Alavi, A. Schwenk (Eds.), {\it Graph Theory, Combinatorics, and Applications}, Vol. 1,
John Wiley \& Sons, Inc., New York, 1995, pp. 311--322.

\bibitem{favaron1996signed}
O. Favaron,
Signed domination in regular graphs,
\DM{158} 
(1-3) (1996), 287--293.

\bibitem{furedi1999signed}
Z. F{\"u}redi, and D. Mubayi,
Signed domination in regular graphs and set systems,
\JCTB{76} 
(2) (1999), 223--239.

\bibitem{henning1999open}
M.A. Henning and P.J. Slater,
Open packing in graphs,
\JCMCC{29} 
(1999), 3--16.


\bibitem{karami2009improved}
H. Karami, A. Khodkar, and S.M. Sheikholeslami,
An improved upper bound for signed edge domination numbers in graphs,
\UM{78} 
(2009), 121--128.

\bibitem{karami2008some}
H. Karami, S. M. Sheikholeslami, and A. Khodkar,
Some notes on signed edge domination in graphs,
\GC{24} 
(1) (2008), 29--35.

\bibitem{matouvsek2000signed}
J. Matou{\v{s}}ek,
On the signed domination in graphs,
\COMBIN{20} 
(1) (2000), 103--108.

\bibitem{xu2001signed}
B. Xu,
On signed edge domination numbers of graphs,
\DM {239} 
(1-3) (2001), 179--189.

\bibitem{xu2005edge}
B. Xu,
On edge domination numbers of graphs,
\DM{294} 
(3) (2005), 311--316.

\bibitem{xu2006two}
B. Xu,
Two classes of edge domination in graphs,
\DAM{154} 
(10) (2006), 1541--1546.

\bibitem{xu2009signed}
B. Xu,
On signed cycle domination in graphs,
\DM{309} 
(4) (2009), 1007--1012.

\bibitem{zelinka1996some}
B. Zelinka,
Some remarks on domination in cubic graphs,
\DM{158} 
(1-3) (1996), 249--255.
\end{thebibliography}
\end{document}